            \def\version{February 20, 2004}           %

\documentclass[reqno,11pt]{amsart}
\usepackage{amsmath, amsthm, a4, latexsym, amssymb}

\makeatletter\@addtoreset{equation}{section}\makeatother
\def\@ref#1{\text{ [#1]}}

\newfam\Bbbfam
\font\tenBbb=msbm10
\font\sevenBbb=msbm7
\font\fiveBbb=msbm5
\textfont\Bbbfam=\tenBbb
\scriptfont\Bbbfam=\sevenBbb
\scriptscriptfont\Bbbfam=\fiveBbb

\def\comment#1{}
\newtheoremstyle{thm}{2ex}{2ex}{\itshape\rmfamily}{}
{\bfseries\rmfamily}{}{1.7ex}{}

\newtheoremstyle{rem}{1.3ex}{1.3ex}{\rmfamily}{}
{\itshape\rmfamily}{}{1.5ex}{}

\newtheorem{theorem}{Theorem}[section]
\newtheorem{lemma}[theorem]{Lemma}
\newtheorem{prop}[theorem] {Proposition}
\newtheorem{cor}[theorem]  {Corollary}

\newtheorem{step}{STEP}
\newtheorem{definition}[theorem] {Definition}

\newcommand{\eq}    {\begin{equation}}
\newcommand{\bel}   {\begin{lemma}}
\newcommand{\el}    {\end{lemma}}
\newcommand{\bep}   {\begin{prop}}
\newcommand{\ep}    {\end{prop}}
\newcommand{\bec}   {\begin{cor}}
\newcommand{\ec}    {\end{cor}}
\newcommand{\bes}   {\begin{step}}
\newcommand{\es}    {\end{step}}
\newcommand{\bea}   {\begin{array}}
\newcommand{\ea}    {\end{array}}

\renewcommand{\subsection}{\secdef \subsct\sbsect}
\newcommand{\subsct}[2][default]{\refstepcounter{subsection}
\nopagebreak\noindent
\vspace{0.5\baselineskip}\noindent
{\flushleft\bf \arabic{section}.\arabic{subsection}~\bf #1  }
\nopagebreak}
\newcommand{\sbsect}[1]{\vspace{0.1cm}\noindent
{\bf #1}\vspace{0.1cm}}

\renewcommand{\subsubsection}{%
\secdef \subsubsect\sbsbsect}
\newcommand{\subsubsect}[2][default]{%
\refstepcounter{subsubsection}
\nopagebreak
\vspace{0.1\baselineskip}
\nopagebreak
{\flushleft
\sffamily\slshape \arabic{section}.\arabic{subsection}.\arabic{subsubsection}
\ %
\sffamily #1\/.}\ }
\newcommand{\sbsbsect}[1]{\vspace{0.1cm}\noindent
{\bf #1}\ }

\newcommand{\eps}       {\varepsilon}

\newcommand{\const}     {{\operatorname {const}\,}}
\newcommand{\esssup}    {{\operatorname {esssup}\,}}
\newcommand{\Prob}      {{\operatorname {Prob}\,}}

\newcommand{\R}     {\mathbb{R}}
\newcommand{\Z}     {\mathbb{Z}}
\newcommand{\N}     {\mathbb{N}}
\renewcommand{\P}   {\mathbb{P}}

\newcommand{\E}     {\mathbb{E}}
\renewcommand{\d}   {\operatorname{d}\!}

\def\1{{\mathchoice {1\mskip-4mu\mathrm l}      
{1\mskip-4mu\mathrm l}
{1\mskip-4.5mu\mathrm l} {1\mskip-5mu\mathrm l}}}

\newcommand{\Ocal}   {{\mathcal O }}

\begin{document}

\title{The parabolic Anderson model}
        
\author[J\"urgen G\"artner and Wolfgang
        K\"onig]{}
\maketitle
\thispagestyle{empty}
\vspace{0.2cm}

\centerline{\sc J\"urgen G\"artner\footnote{Technische Universit\"at Berlin, Institut f\"ur Mathematik,~MA7-5,
Stra\ss e des 17.~Juni~136, D-10623 Berlin, Germany.}$^,$\footnote{Partially supported by the 
German Science Foundation, Schwerpunkt project SPP 1033}} 
\centerline{\sc Wolfgang  K\"onig\footnotemark[1]$^,$\footnotemark[2]$^,$\footnote{{\tt koenig@math.tu-berlin.de}}} 

\vspace{0.8cm}

\centerline{\small(\version)}\vspace{6ex}

\begin{quote}{\small {\bf Abstract:} This is a survey on the intermittent 
behavior of the parabolic \mbox{Anderson} model, which is the Cauchy problem 
for the heat equation with random potential on the lattice $\Z^d$. We first 
introduce the model and give heuristic explanations of the long-time behavior 
of the solution, both in the annealed and the quenched setting for 
time-independent potentials. We thereby consider examples of 
potentials studied in the literature. In the particularly 
important case of an i.i.d.\ potential with double-exponential tails 
we formulate the asymptotic results in detail. Furthermore, 
we explain that, under mild regularity assumptions, there are only four 
different universality classes of asymptotic behaviors. Finally, we study 
the moment Lyapunov exponents for space-time homogeneous catalytic potentials 
generated by a \mbox{Poisson} field of random walks.
}
\end{quote}


\vfill



\pagebreak

\setcounter{section}{0}
\section{Introduction and heuristics}

\subsection{Evolution of spatially distributed systems in random media}

\noindent
One of the often adequate and frequently used methods for studying 
the evolution of spatially distributed systems under the influence of a 
random medium 
is {\em homogenization}. After rescaling, the system, modeled by 
partial differential equations with random coefficients, is approached by 
a system with `properly averaged' deterministic coefficients, see e.g.\ 
\cite{ZKO94}. But there are simple and important situations when random 
systems exhibit effects which cannot be recovered by such deterministic 
approximations and related fluctuation corrections. This concerns, 
in particular, {\em localization\/} effects for non-reversible 
random walks in random environment \cite{S82} and for the electron 
transport in disordered media \cite{And58}. 

Another such effect is that of {\em intermittency}. Roughly speaking, 
intermittency means that the solution of the system develops pronounced 
spatial structures on islands located far from each other that,
in one or another sense, deliver the main output to the system. One of the 
sources  of interest is magnetohydrodynamics and, in particular, the 
investigation of the induction equation with incompressible random velocity 
fields \cite{Z84}, \cite{ZMRS87}. Another source are simple mathematical 
models such as the random \mbox{Fisher}-\mbox{Eigen} equation that have been 
used to derive caricatures of \mbox{Darwinian} evolution principles 
\cite{EEEF84}. 

One of the simplest and most basic models exhibiting the effect of 
intermittency is the Cauchy problem for the spatially discrete heat equation 
with a random potential: 
\begin{equation}\label{Anderson0}
\begin{array}{rcll}
\displaystyle
\partial_t \,u(t,x)\!\!\! &=&\!\!\!\kappa \Delta u(t,x)+\xi(t,x) 
u(t,x),\qquad &(t,x)\in(0,
\infty)\times {\mathbb Z}^d,\\
u(0,x)\!\!\!&=&\!\!\!\!u_0(x),&x\in{\mathbb Z}^d.
\end{array}
\end{equation}
Here $\kappa>0$ is a diffusion constant, $\Delta$ denotes the discrete 
\mbox{Laplacian}, 
$$
\Delta f(x)=\sum_{y\colon |y-x|=1}[f(y)-f(x)],
$$
$\xi$ is a space-time homogeneous ergodic random potential, and $u_0$ is a 
nonnegative initial function. Problem (\ref{Anderson0}) is often called 
{\em parabolic problem for the \mbox{Anderson} model\/} or 
{\em parabolic \mbox{Anderson} model\/} 
(abbreviated {\em PAM\/}). As simplest localized initial datum one may 
take $u_0=\delta_0$, and as non-local initial datum $u_0=\1$. 
In the latter case, the solution $u(t,\cdot)$ is spatially homogeneous and
ergodic for each $t$. Let us remark that the solution $u(t,x)$ of 
(\ref{Anderson0}) allows the interpretation as average number of 
particles at site $x$ at time $t$ for branching random walks in random 
media given a realization of the medium $\xi$, cf.\ \cite{CM94} and 
the remarks in the next subsection.

\subsection{The PAM with time-independent potential}\label{Pamintro}

\noindent
In the particular case when the potential $\xi(t,x)=\xi(x)$ is 
time-independent, the large-time behavior of the solution $u$ to the 
PAM (\ref{Anderson0}) is determined by the spectral properties of the 
\mbox{Anderson} \mbox{Hamiltonian}
\begin{equation}\label{AndHam}
\mathcal H = \kappa\Delta + \xi
\end{equation}
and therefore closely related to the mentioned localization of the 
electron transport. Namely, since (under natural assumptions on $\xi$) 
the upper part of the spectrum of 
$\mathcal H$ in $\ell^2(\mathbb Z^d)$ is a pure point spectrum 
\cite{FMSS85}, \cite{AM93}, the solution $u$ admits the spectral 
representation
\begin{equation}\label{SR}
u(t,\cdot) = \sum_i e^{\lambda_i t}\left(v_i,u_0\right) v_i(\cdot)
\end{equation}
with respect to the random eigenvalues $\lambda_i$ and the corresponding 
exponentially localized random eigenfunctions $v_i$. 
(For simplicity we ignore the possible occurrence of a continuous central part 
of the spectrum.) As $t$ increases 
unboundedly, only summands with larger and larger eigenvalues will 
contribute to (\ref{SR}), and the corresponding eigenfunctions are 
expected to be localized more and more far from each other. Hence, for 
large $t$, the solution $u(t,\cdot)$ indeed looks like a weighted 
superposition of high peaks concentrated on distant islands. 

A mathematically rigorous understanding of the nature of the spectrum 
of the \mbox{Anderson} \mbox{Hamiltonian} and the random 
\mbox{Schr\"odinger} operator is still far from being complete. For 
an overview about some recent developments we refer to the surveys 
\cite{BKS04} and \cite{LMW04} in this proceedings volume. 
The spectral results 
obtained so far do not yet seem directly applicable to answer the 
crucial questions about intermittency. A direct spectral approach 
clearly fails for space-time dependent potentials. 

In this survey we present a part of the results about intermittency for the 
PAM which have been obtained by use of more intrinsic probabilistic methods. 
In the next subsections we stick to the PAM with time-independent 
potential and localized initial datum:
\begin{equation}\label{Anderson4}
\begin{array}{rcll}
\displaystyle
\partial_t \,u(t,x)\!\!\! &=&\!\!\!\kappa \Delta u(t,x)+\xi(x) 
u(t,x),\qquad &(t,x)\in(0,
\infty)\times {\mathbb Z}^d,\\
u(0,x)\!\!\!&=&\!\!\!\delta_0(x),&x\in{\mathbb Z}^d.
\end{array}
\end{equation}
We assume throughout that $\xi=(\xi(x))_{x\in{\mathbb Z}^d}$ is a field 
of i.i.d.\ random variables with finite positive exponential moments. 
Under these basic assumptions, $u(t,x)$ has moments of all orders. 
 
The solution $u$ to (\ref{Anderson4}) describes a random particle flow in 
$\Z^d$ in the presence 
of  random sources (lattice sites $x$ with $\xi(x)>0$) and random 
sinks (sites $x$ with $\xi(x)<0$).\footnote{Sites $x$ with 
$\xi(x)=-\infty$ may be allowed and interpreted as (`hard') traps or 
obstacles, sites with 
$\xi(x)\in(-\infty,0)$ are sometimes called `soft' traps.} 
Two competing effects are present: the diffusion mechanism governed 
by the \mbox{Laplacian}, and the local growth governed by the 
potential. The diffusion tends to make the random field
$u(t,\cdot)$ flat, whereas the random potential $\xi$ has a tendency to make 
it irregular.

The solution $u$ to (\ref{Anderson4}) also admits a branching particle 
dynamics interpretation. 
Imagine that initially, at time $t=0$, there is a single particle at the 
origin, and all other sites are vacant. This particle moves according to a 
continuous-time symmetric random walk with generator $\kappa\Delta$. 
When present at site $x$, the particle is split into two particles with rate $\xi_+(x)$ and is killed with rate $\xi_-(x)$, where $\xi_+=(\xi_+(x))_{x\in\Z^d}$ and $\xi_-=(\xi_-(x))_{x\in\Z^d}$ are independent random i.i.d.\ fields 
($\xi_-(x)$ may attain the value $\infty$). Every particle continues from its birth site in the same way as the parent particle, and their movements are independent. Put $\xi(x)=\xi_+(x)-\xi_-(x)$. Then, given $\xi_-$ and $\xi_+$, 
the expected number of particles present at the site $x$ at 
time $t$ is equal to $u(t,x)$. Here the expectation is taken over the particle motion and over the splitting resp.~killing mechanism, but not over the random medium $(\xi_-,\xi_+)$. 

A very useful standard tool for the probabilistic investigation  of 
(\ref{Anderson4}) is the well-known {\it Feynman-Kac formula\/}
for the solution $u$, which (after time-reversal) reads 
\begin{equation}\label{FKform}
u(t,x)=\E_0 \Bigl[\exp\Bigl\{\int_0^t\xi(X(s))\, ds\Bigr\}\delta_x(X(t))\Bigr],\qquad (t,x)\in[0,\infty)\times\mathbb Z^d,
\end{equation}
where $(X(s))_{s\in [0,\infty)}$ is continuous-time random walk
on ${\mathbb Z}^d$ with generator $\kappa\Delta$ starting at $x\in {\mathbb Z}^d$ under $\E_x$.  

Our main interest concerns the large-time behavior of the random field $u(t,\cdot)$. In particular, we consider the total mass, i.e., the random variable
\begin{equation}\label{Udef}
U(t)=\sum_{x\in\Z^d}u(t,x)=\E_0 \Bigl[\exp\Bigl\{\int_0^t\xi(X(s))\, ds\Bigr\}\Bigr],\qquad t>0. 
\end{equation}
Note that $U(t)$ coincides with the value $\hat{u}(t,0)$ of the solution 
$\hat{u}$ to the parabolic equation \eqref{Anderson4} with initial 
datum $\hat{u}(0,\cdot)=\1$ instead of $u(0,\cdot)=\delta_0$. One 
should have in mind that, because of this, our considerations below also 
concern the large-time asymptotics of $\hat{u}$. 

We ask the following questions:
\begin{enumerate}
\item What is the asymptotic behavior of $U(t)$ as $t\to\infty$?

\item Where does the main mass of $u(t,\cdot)$ stem from? What are the regions that contribute most to $U(t)$? What are these regions determined by? How many of them are there and how far away are they from each other?

\item What do the typical shapes of the potential $\xi(\cdot)$ and of the solution $u(t,\cdot)$ look like in these regions?
\end{enumerate}

We call the regions that contribute the overwhelming part to the total mass $U(t)$ {\em relevant islands\/} or {\em relevant regions}. The notion of {\em intermittency\/} states that there does exist a small number of relevant islands which are far away from each other and carry asymptotically almost all the total  mass $U(t)$ of $u(t,\cdot)$. See Section~\ref{IntSect} for details. 

This effect may also be studied from the point of view of {\em typical 
paths\/} $X(s)$, $s\in[0,t]$, giving the main contribution to the expectation 
in the \mbox{Feynman}-\mbox{Kac} formula (\ref{Udef}).  On the one hand, 
the random walker $X$ should move quickly and as far as possible through 
the potential landscape to reach a region of exceptionally 
high potential and then stay there up to time $t$. This will make the integral 
in the exponent on the right of (\ref{Udef}) large. On the other hand, 
the probability to reach such a distant potential peak up to $t$ may be rather 
small. Hence, the first order contribution to $U(t)$ comes from paths that 
find a good compromise between the high potential values and the far distance. 
This contribution is given by the height of the peak. The second order 
contribution to $U(t)$ is determined by the precise manner in which the 
optimal walker moves within the potential peak, and this depends on the 
geometric properties of the potential in that peak. 

It is part of our study to understand the effect of intermittency for the 
parabolic Anderson model in great detail.
We distinguish between the so-called {\em quenched\/} setting, where we consider $u(t,\cdot)$ 
almost surely with respect to the medium $\xi$, and the {\em annealed\/} one, where we average with respect to $\xi$. It is clear that the quantitative 
details of the answers to the above questions strongly depend on the distribution of the field $\xi$ (more precisely, on the upper tail of the distribution of the random variable $\xi(0)$), and that different phenomena occur in the quenched and the annealed settings. 

It will turn out that there is a universal picture present in the asymptotics 
of the parabolic Anderson model. Inside the relevant islands, after appropriate 
vertical shifting and spatial rescaling, the potential $\xi$ will turn out to 
asymptotically approximate a universal, non-random shape, $V$, which is determined 
by a characteristic variational problem. The absolute height of the potential 
peaks and the diameter of the relevant islands are asymptotically determined by 
the upper tails of the random variable $\xi(0)$, while the number of the islands 
and their locations are random. Furthermore, after multiplication with an 
appropriate factor and rescaling, also the solution $u(t,\cdot)$ approaches a 
universal shape on these islands, namely the principal eigenfunction of 
the \mbox{Hamiltonian} $\kappa\Delta+V$  with $V$ the above universal potential 
shape. Remarkably, there are only four universal classes of potential shapes 
for the PAM in (\ref{Anderson4}), see Section~\ref{Universality} for details. 

For a general discussion we refer to the monograph by 
Carmona and Molchanov~\cite{CM94}, the lectures by
Molchanov~\cite{M94}, and also to the results by Sznitman about the 
important (spatially continuous) 
case of bounded from above \mbox{Poisson}-like potentials summarized in his
monograph \cite{S98}. A discussion from a physicist's and a chemist's point of view in the particular case
of trapping problems (see also Section~\ref{survival4} below), 
including a survey on related mathematical models and a 
collection of open problems, is provided in \cite{HW94}.
A general mathematical background for the PAM is provided in \cite{GM90}.

\subsection{Intermittency}\label{IntSect}

\noindent
As before, let $\hat{u}$ denote the solution to the equation in 
(\ref{Anderson4}) with initial datum $\hat{u}(0,\cdot)=\1$, 
but now with a homogeneous ergodic potential $\xi=(\xi(x))_{x\in\mathbb Z^d}$. 
Assume that all positive exponential moments of $\xi(0)$ are finite.   
Let $\Prob(\cdot)$ and $\langle\cdot\rangle$ denote probability and 
expectation w.r.t.\ $\xi$. 

A first, rough, mathematical approach to intermittency consists in a 
comparision of the growth of subsequent moments of the ergodic field 
$\hat{u}(t,\cdot)$ as $t\to\infty$. Define
$$
\Lambda_p(t) = \log \left\langle\hat{u}(t,0)^p\right\rangle,
\qquad p\in\mathbb N, 
$$  
and write $f\ll g$ if $\lim_{t\to\infty}[g(t)-f(t)]=\infty$. 

\begin{definition}\label{Defintermit}\rm
For $p\in\mathbb N\setminus\{1\}$, the homogeneous ergodic field 
$\hat{u}(t,\cdot)$ is called {\em $p$-intermittent\/} as 
$t\to\infty$, if
\begin{equation}\label{intermitdef}
\frac{\Lambda_{p-1}}{p-1} \ll \frac{\Lambda_p}{p}.
\end{equation}
\end{definition}

\noindent
Note that, by \mbox{H\"older's} inequality, always 
$\Lambda_{p-1}/(p-1)\le\Lambda_p/p$. 
If the finite moment \mbox{Lyapunov} exponents 
$$
\lambda_p = \lim_{t\to\infty} \frac{1}{t}\Lambda_p(t),
\qquad p\in\mathbb N, 
$$
exist, then the strict inequality 
$\lambda_{p-1}/(p-1)<\lambda_p/p$ implies $p$-intermittency. Such a 
comparison of the moment \mbox{Lyapunov} exponents has first been used in the 
physics literature to study intermittency, cf.\ 
\cite{ZMRS87}, \cite{ZMRS88}. We will use this approach in 
Section~\ref{Catsection}.

To explain the meaning of Definition~\ref{Defintermit}, assume 
(\ref{intermitdef}) for some $p\in\mathbb N\setminus\{1\}$ 
and choose a level function $\ell_p$ such that 
$\Lambda_{p-1}/(p-1)\ll \ell_p \ll \Lambda_p/p$. Then, by \mbox{Chebyshev's} 
inequality,
$$
\Prob\!\!\left(\hat{u}(t,0)>e^{\ell_p(t)}\right) \le 
e^{-(p-1)\ell_p(t)} \left\langle\hat{u}(t,0)^{p-1}\right\rangle
= \exp\!\left\{\Lambda_{p-1}(t)-(p-1)\ell_p(t)\right\},
$$
and the expression on the right converges to zero as $t\to\infty$. 
In other words, the density of the homogeneous point process 
$$
\Gamma(t) = \left\{x\in\mathbb Z^d\colon \hat{u}(t,x)>e^{\ell_p(t)}\right\}
$$
vanishes asymptotically as $t\to\infty$. On the other hand,
\begin{eqnarray*}
\left\langle \hat{u}(t,0)^p \,\1\!\left\{\hat{u}(t,0) \le e^{\ell_p(t)}
\right\}\right\rangle \le e^{p\ell_p(t)}
= e^{p\ell_p(t)-\Lambda_p(t)}\left\langle\hat{u}(t,0)^p\right\rangle
= o\left(\left\langle\hat{u}(t,0)^p\right\rangle\right)
\end{eqnarray*}
and, consequently,
$$
\left\langle\hat{u}(t,0)^p\right\rangle
\,\sim\, \left\langle\hat{u}(t,0)^p\,\1\!\left\{\hat{u}(t,0)>e^{\ell_p(t)}
\right\}\right\rangle
$$
as $t\to\infty$. Hence, by \mbox{Birkhoff's} ergodic theorem, for large $t$ 
and large centered boxes $B$ in $\mathbb Z^d$, 
$$
|B|^{-1} \sum_{x\in B}\hat{u}(t,x)^p 
\,\approx\, |B|^{-1}\!\sum_{x\in B\cap \Gamma(t)} \hat{u}(t,x)^p.
$$
This means that the $p$-th moment $\langle\hat{u}(t,0)^p\rangle$ is 
`generated' by the high peaks of $\hat{u}(t,\cdot)$ on the `thin' set 
$\Gamma(t)$ and therefore  indicates the presence of intermittency 
in the above verbal sense. Unfortunately, this approach does not reflect 
the geometric structure of the set $\Gamma(t)$. This set might consist of 
islands or, e.g., have a net-like structure. 

\begin{theorem}
If $\xi=(\xi(x))_{x\in\mathbb Z^d}$ is a non-deterministic field of i.i.d.\ 
random variables with $\langle e^{t\xi(0)}\rangle<\infty$ for all $t>0$, 
then the solution $\hat{u}(t,\cdot)$ is $p$-intermittent for all 
$p\in\mathbb N\setminus\{1\}$. 
\end{theorem}

\noindent
This is part of Theorem~3.2 in \cite{GM90}, where,
for general homogeneous ergodic potentials $\xi$, necessary and sufficient 
conditions for $p$-intermittency of $\hat{u}(t,\cdot)$ have been given 
in spectral terms of the \mbox{Hamiltonian} (\ref{AndHam}).

\subsection{Annealed second order asymptotics}\label{heur4}

\noindent 
Let us discuss, on a heuristic level, what the asymptotics of the 
moments of $U(t)$ are determined by, and how they can be described. 
For simplicity we restrict ourselves to the first moment.

The basic observation is that, as a consequence of the spectral 
representation (\ref{SR}), 
\begin{equation}\label{appreigenv}
U(t) \approx e^{t\lambda_t(\xi)}
\end{equation}
(in the sense of logarithmic equivalence), where $\lambda_t(\varphi)$ denotes 
the principal (i.e., largest) eigenvalue 
of the operator $\kappa\Delta+\varphi$ with zero boundary condition in the 
`macrobox' 
$B_t=[-t,t]^d\cap\Z^d$. Hence, we have to understand the large-time 
behavior of the exponential moments of the principal eigenvalue of the 
\mbox{Anderson} \mbox{Hamiltonian} $\mathcal H$ in a large, time-dependent 
box.

It turns out that the main contribution to 
$\langle e^{t\lambda_t(\xi)}\rangle$ comes from realizations of the 
potential $\xi$ having high peaks on distant islands of some radius of order 
$\alpha(t)$  
that is much smaller than $t$. But this implies that  $\lambda_t(\xi)$ is 
close to the principal eigenvalue of $\mathcal H$ on one of these islands. 
Therefore, since the number of subboxes of $B_t$ of radius of order 
$\alpha(t)$ grows only polynomial in $t$ and $\xi$ is spatially homogeneous, 
we may expect that 
$$
\left\langle e^{t\lambda_t(\xi)} \right\rangle 
\approx \left\langle e^{t\lambda_{R\alpha(t)}(\xi)} \right\rangle 
$$
for $R$ large as $t\to\infty$. 
 
The choice of the scale function $\alpha(t)$ depends on asymptotic 
`stiffness' properties of the potential, more precisely of its tails at its 
essential supremum, and is determined by a large deviation principle, 
see \eqref{LDPpot} below. 
In Section~\ref{examples} we shall see examples of potentials 
such that $\alpha(t)$ tends to 0, to $\infty$, or stays bounded and bounded 
away from zero as $t\to\infty$. In the present heuristics, we shall assume 
that $\alpha(t)\to\infty$, which implies the necessity of a spatial 
rescaling. In particular, after rescaling, the main quantities and objects 
will be described in terms of the continuous counterparts of the discrete 
objects we started with, i.e., instead of the discrete Laplacian, 
the continuous Laplace operator appears etc. The following heuristics can 
also be read in the case where $\alpha(t)\equiv 1$ by keeping the discrete 
versions for the limiting objects.

The optimal behavior of the field $\xi$ in the `microbox' $B_{R\alpha(t)}$ 
is to approximate a certain (deterministic) shape $\varphi$ after appropriate 
spatial scaling and vertical shifting. 
It easily follows from the Feynman-Kac formula \eqref{Udef} that  
$$
e^{H(t)-2d\kappa t} \le \langle U(t) \rangle \le e^{H(t)}, 
$$
where
\begin{equation}\label{Hdef}
H(t)=\log \bigl\langle e^{t\xi(0)}\bigr\rangle,\qquad t>0,
\end{equation}
denotes the {\em cumulant generating function\/} of $\xi(0)$ 
(often called {\em logarithmic moment generating 
function\/}). Hence, the peaks of $\xi(\cdot)$ mainly contributing to 
$\langle U(t)\rangle$ have height of order $H(t)/t$. Together with 
\mbox{Brownian} scaling this leads to the ansatz
\begin{equation}
\overline\xi_t(\cdot)=\alpha(t)^2\Bigl[\xi\bigl(\lfloor \,\cdot\,\alpha(t)
\rfloor\bigr) -\frac {H(t)}t\Bigr],
\end{equation}
for the spatially rescaled and vertically shifted potential in the 
cube $Q_R=(-R,R)^d$. Now the idea is that the main contribution to 
$\langle U(t)\rangle$ comes from fields that are shaped in such a way that 
$\overline\xi_t\approx \varphi$ in $Q_R$, for some $\varphi\colon Q_R\to\R$, 
which has to be chosen optimally. Observe that
\begin{equation}\label{xiscale}
\overline\xi_t\approx \varphi\quad \text{in }Q_R\qquad\Longleftrightarrow
\qquad \xi(\cdot)\approx \frac {H(t)}t+{\textstyle{\frac 1{\alpha(t)^2}
\varphi\bigl(\frac{\cdot}{\alpha(t)}\bigr)}} \quad\text{in }B_{R\alpha(t)}.
\end{equation}
Let us calculate the contribution to $\langle U(t)\rangle$ coming from such 
fields. 
Using \eqref{appreigenv}, we obtain
\begin{equation}
\left\langle U(t)\,\1\{\overline\xi_t\approx \varphi\text{ in } Q_R\}
\right\rangle\approx e^{H(t)} \exp\Bigl\{t\lambda_{R\alpha(t)}
\bigl({\textstyle{\frac 1{\alpha(t)^2}\varphi\bigl(\frac{\cdot}{\alpha(t)}
\bigr)}}\bigr)\Bigr\}\,\Prob\bigl(\overline\xi_t\approx \varphi\text{ in } 
Q_R\bigr).
\end{equation}
The asymptotic scaling properties of the discrete Laplacian, $\Delta$, 
imply that
\begin{equation}
\lambda_{R\alpha(t)}\bigl({\textstyle{\frac 1{\alpha(t)^2}
\varphi\bigl(\frac{\cdot}{\alpha(t)}\bigr)}}\bigr)\approx 
\frac{1}{\alpha(t)^2} \lambda^{\rm c}_R(\varphi),
\end{equation}
where $\lambda^{\rm c}_R(\varphi)$ denotes the principal eigenvalue of 
$\kappa\Delta^{\rm c}+\varphi$ in the cube $Q_R$ with zero boundary condition, 
and $\Delta^{\rm c}$ is the usual `continuous' Laplacian. This leads to 
\begin{equation}\label{appr3}
\left\langle U(t)\,\1\{\overline\xi_t\approx \varphi\text{ in } Q_R\}
\right\rangle\approx e^{H(t)} 
\exp\Bigl\{\frac t{\alpha(t)^2}\lambda^{\rm c}_R(\varphi)\Bigr\}\,
\Prob\bigl(\overline\xi_t\approx \varphi\text{ in } Q_R\bigr).
\end{equation}
In order to achieve a balance between the second and the third factor on the 
right, it is necessary that the logarithmic decay rate of the considered 
probability is $t/\alpha(t)^2$. One expects to have a large deviation 
principle for the shifted, rescaled field, which reads
\begin{equation}\label{LDPpot}
\Prob\bigl(\overline\xi_t\approx \varphi\text{ in } Q_R\bigr)
\approx \exp\Bigl\{-\frac t{\alpha(t)^2} I_R(\varphi)\Bigr\},
\end{equation}
where the scale $\alpha(t)$ has to be determined in such a way that the 
rate function $I_R$ is non-degenerate. Now substitute \eqref{LDPpot} into 
\eqref{appr3}. Then the \mbox{Laplace} method tells us that the exponential 
asymptotics of $\langle U(t)\rangle$ is equal to the one of 
$\bigl\langle U(t)\,\1\{\overline\xi_t\approx \varphi\text{ in } Q_R\}$ 
with optimal $\varphi$. Hence, optimizing on $\varphi$ and remembering that 
$R$ is large, we arrive at
\begin{equation}\label{appr4}
\bigl\langle U(t)\bigr\rangle\approx e^{H(t)} 
\exp\Bigl\{-\frac t{\alpha(t)^2}\chi\Bigr\},
\end{equation}
where the constant $\chi$ is given in terms of the characteristic 
variational problem
\begin{equation}\label{chideffirst}
\chi=\lim_{R\to\infty}\inf_{\varphi\colon Q_R\to\R}\bigl[I_R(\varphi)
-\lambda^{\rm c}_R(\varphi)\bigr].
\end{equation}
The first term on the right of \eqref{appr4} is determined by the absolute 
height of the typical realizations of the potential and the second contains 
information about the shape of the potential close to its maximum in spectral 
terms of 
the \mbox{Anderson} \mbox{Hamiltonian} $\mathcal H$ in this region. 
More precisely, those realizations of $\xi$ with 
$\overline\xi_t\approx \varphi_*\text{ in } Q_R$ for large $R$ and 
$\varphi_*$ a minimizer in the variational formula in \eqref{chideffirst} 
contribute most to $\langle U(t)\rangle$. In particular, the geometry of the 
relevant potential peaks is hidden via $\chi$ in the second asymptotic term 
of $\langle U(t)\rangle$.

\subsection{Quenched second order asymptotics}\label{Quenchedheu}

\noindent 
Here we explain, again on a heuristic level, the almost sure asymptotics of 
$U(t)$ as $t\to\infty$. Because of \eqref{appreigenv}, it suffices to study  
the asymptotics of the principal eigenvalue $\lambda_t(\xi)$. 

Like for the annealed asymptotics, the main contribution to 
$\lambda_t(\xi)$ comes from islands whose radius is of a certain 
deterministic, time-depending order, which we denote $\widetilde\alpha(t)$. 
As $t\to\infty$, the scale function $\widetilde\alpha(t)$ tends to zero, one, 
or $\infty$, respectively, if the scale function $\alpha(t)$ for the moments 
tends to these respective values (see also \eqref{alphadef} below). 
However, $\widetilde\alpha(t)$ is roughly of logarithmic order in $\alpha(t)$ 
if $\alpha(t)\to\infty$, hence it is {\em much\/} smaller than $\alpha(t)$.  

The relevant islands (`microboxes') have radius $R\widetilde\alpha(t)$, where 
$R$ is 
chosen large. Let $z\in B_t$ denote the (certainly random) center of one of 
these islands $\widetilde B=z+B_{R\widetilde\alpha(t)}$ meeting 
the two requirements (1) the potential $\xi$ is very large in $\widetilde B$ 
and (2) $\xi$ has an optimal shape within $\widetilde B$. This is further 
explained as follows. Let $h_t=\max_{B_t} \xi$ be the maximal potential value 
in the large box $B_t$. (Then $h_t$ is  a priori random, but well approximated 
by deterministic asymptotics, which can be deduced from asymptotics of 
$H(t)$.) Then $\xi-h_t$ is roughly of finite order within the relevant 
`microbox' $\widetilde B$. Furthermore, $\xi-h_t$ should approximate a fixed 
deterministic shape in $\widetilde B$. Hence, we consider the shifted and 
rescaled field in the box $\widetilde B$, 
\begin{equation}\label{shiftresc}
\overline \xi_{t}(\cdot)=\widetilde\alpha(t)^2\Bigl[\xi\bigl(z+\cdot\,\widetilde\alpha(t)\bigr)-h_t\Bigr],\qquad\mbox{in }Q_R=(-R,R)^d.
\end{equation}
Note that 
\begin{equation}\label{potapr}
\overline \xi_{t}\approx \varphi\quad\mbox{ in }Q_R\qquad\Longleftrightarrow\qquad \xi(z+\cdot)\approx h_t+{\textstyle{\frac{1}{\widetilde\alpha(t)^2}\varphi\bigl(\frac\cdot{\widetilde\alpha(t)}\bigr)}}\quad\mbox{ in }\widetilde B-z.
\end{equation}

A crucial Borel-Cantelli argument shows that, for a given shape $\varphi$, 
with probability one, for any $t$ sufficiently large, there does exist 
at least one box $\widetilde B$ having radius $R\widetilde\alpha(t)$ such that 
the event $\{\overline \xi_{t}\approx \varphi\mbox{ in }Q_R\}$ occurs
if $I_R(\varphi)<1$, where $I_R$ is the rate function of the 
large deviation principle in \eqref{LDPpot}. If $I_R(\varphi)>1$, then this 
happens with probability 0. For the Borel-Cantelli argument 
to work, one needs the scale function $\widetilde\alpha(t)$ to be  
defined in terms of the annealed scale function $\alpha(t)$ in the following 
way:
\begin{equation}\label{alphadef}
\frac{\widetilde \alpha(t)}{\alpha(\widetilde\alpha(t))^2}=d\log t,
\end{equation}
i.e., $\widetilde\alpha(t)$ is the inverse of the map $t\mapsto t/\alpha(t)^2$, evaluated at $d\log t$. Note that the growth of $\widetilde\alpha(t)$ is roughly of logarithmic order of the growth of $\alpha(t$), i.e., if the annealed relevant islands grow unboundedly, then the quenched relevant islands also grow unboundedly, but with much smaller velocity.

Hence, with probability one, for all large $t$, there is at least one box 
$\widetilde B$ in which the potential looks like the function on the right 
of \eqref{potapr}. The contribution to $\lambda_t(\xi)$ coming from one of 
the boxes $\widetilde B$ is equal to the associated principal eigenvalue 
\begin{equation}\label{lambdascal}
\lambda_{\widetilde B-z}\bigl(h_t+  
{\textstyle{\frac{1}{\widetilde\alpha(t)^2}
\varphi\bigl(\frac\cdot{\widetilde\alpha(t)}\bigr)}}\bigr)
\approx h_t+ \frac 1{\widetilde\alpha(t)^2}\lambda^{\rm c}_R(\varphi),
\end{equation}
where we recall that $\lambda^{\rm c}_R(\varphi)$ is the principal 
\mbox{Dirichlet} eigenvalue 
of the operator $\kappa\Delta^c+\varphi$ in the `continuous' cube $Q_R$. 
Obviously, $\lambda_t(\xi)$ is asymptotically not smaller than the expression 
on the right of \eqref{lambdascal}. In terms of the Feynman-Kac formula in 
\eqref{FKform}, this lower estimate is obtained by inserting the indicator on 
the event that the random path moves quickly to the box $\widetilde B$ and 
stays all the time until $t$ in that box. 

It is an important technical issue to show that, asymptotically as 
$t\to\infty$, $\lambda_t(\xi)$ is also estimated from {\em above\/} by the 
right hand side of \eqref{lambdascal}, if $\varphi$ is optimally chosen, i.e., 
if $\lambda^{\rm c}_R(\varphi)$ is optimized over all admissible $\varphi$ 
and on $R$. This implies that the almost sure asymptotics of $U(t)$ are given 
as 
\begin{equation}
\frac {1}{t}\log U(t)\approx \lambda_t(\xi)
\approx h_t -\frac{1}{\widetilde\alpha(t)^{2}}\widetilde\chi,
\qquad t\to\infty,
\end{equation}
where $\widetilde \chi$ is given in terms of the characteristic variational 
problem
\begin{equation}\label{chitilde}
\widetilde\chi=\lim_{R\to\infty}\inf_{\varphi\colon Q_R\to\R, 
I_R(\varphi)<1} \left[-\lambda^{\rm c}_R(\varphi)\right].
\end{equation}
This ends the heuristic derivation of the almost sure asymptotics of $U(t)$. 
Like in the annealed case, there are two terms, which describe the absolute 
height of the potential in the `macrobox' $B_t$, and the shape of the 
potential in the relevant `microbox' $\widetilde B$, more precisely the 
spectral properties of $\kappa\Delta+\xi$ in that microbox. The interpretation 
is that, for $R$ large  and $\varphi_*$ an approximate minimizer in 
\eqref{chitilde}, the main contribution to $U(t)$ comes from a small box 
$\widetilde B$ in $B_t$, with radius $R\widetilde\alpha(t)$, in which the 
shifted and rescaled potential $\overline \xi_t$ looks like $\varphi_*$. 
The condition $I_R(\varphi_*)<1$ guarantees the existence of such a box, 
and $\lambda_R(\varphi_*)$ quantifies the contribution from that box.

Let us remark that the variational formulas in \eqref{chitilde} and 
\eqref{chideffirst} are in close connection to each other. In particular, 
it can be shown that the minimizers of \eqref{chitilde} are rescaled versions 
of the minimizers of \eqref{chideffirst}. This means that, up to rescaling, 
the optimal potential shapes in the annealed and in the quenched setting are 
identical.

\subsection{Geometric picture of intermittency}\label{Interheu}

\noindent 
In this section we explain the geometric picture of intermittency, 
still on a heuristic level.

The heuristics for the total mass of $u(t,\cdot)$ in Section~\ref{Quenchedheu}
makes use of only {\it one\/} of the relevant islands 
$\widetilde B$ in which the potential is optimally valued and shaped. 
In order to describe the entire function $u(t,\cdot)$, one has to take into 
account a certain (random) number of such islands. Let $n(t)$ denote their 
number, and let $z_1,z_2,\dots,z_{n(t)}\in B_t$ denote the centers of these 
relevant microboxes $B_1,B_2,\dots,B_{n(t)}$, whose radii are equal to 
$R\widetilde\alpha(t)$. Then, almost surely, 
\begin{equation}
U(t)=\sum_{x\in\Z^d}u(t,x)\approx \sum_{i=1}^{n(t)}\sum_{x\in B_i}u(t,x),
\qquad \text{as $t\to\infty$,}
\end{equation}
i.e., asymptotically the total mass of the random field $u(t,\cdot)$ stems 
only from the unions of the relevant islands, $B_1,\dots,B_{n(t)}$. These 
islands are far away from each other. On each of them, the shifted and 
rescaled potential $\overline\xi_t$, see \eqref{shiftresc}), looks 
approximately like a minimizer $\varphi_*$ of the variational problem in 
\eqref{chitilde}. In particular, it has an asymptotically 
{\em deterministic\/} shape. This is the universality in the potential 
landscape: the height and the (appropriately rescaled) shape of the 
potential on the relevant islands are deterministic, but their location and
number are random.  

The shape of the solution, $u(t,\cdot)$, on each of  the relevant islands 
also approaches a universal deterministic shape, namely a time-dependent multiple 
of the principal eigenfunction of the operator $\kappa\Delta+\varphi_*$.

\section{Examples of potentials}\label{examples}

\subsection{Double-exponential distributions} \label{doubleexp}

\noindent Consider a distribution which lies in the vicinity of the 
{\em double-exponential distribution}, i.e.,
\begin{equation}\label{doubleexp4}
\Prob(\xi(0)>r)\approx\exp\bigl\{- e^{r/\varrho}\bigr\},
\qquad r\to\infty,
\end{equation}
with $\varrho\in(0,\infty)$ a parameter. It turns out \cite{GM98} 
that this class of potentials constitutes a critical class in the sense 
that the radius of the relevant islands stays finite as $t\to\infty$. 
This is related to the characteristic property of the double-exponential 
distribution that
$$
\Prob(\xi(x)>h)\approx \Prob(\xi(y)>h-\varrho\log 2,\,\xi(z)>h-\varrho\log 2),
$$
meaning that single-site potential peaks of height $h\gg 1$ occur with the 
same frequency as two-site potential peaks with height of the same order. 
Hence, no spatial rescaling is necessary, and we put $\alpha(t)=1$. 
In Sections~\ref{Annres}--\ref{Interres} below we shall describe our results 
for this type of potentials more closely.

For the boundary cases $\varrho=\infty$ and $\varrho=0$ (`beyond' and 
`on this side of' the double-exponential distribution, respectively), 
\cite{GM98} argued that the boundary cases $\alpha(t)\downarrow 0$ and 
$\alpha(t)\to\infty$ occur. In other words, the fields beyond the 
double-exponential (which includes, e.g., Gaussian fields) are simple in the 
sense that the main contribution comes from islands consisting of single 
lattice sites. Unbounded fields that are in the vicinity of the case
$\varrho=0$ are called `almost bounded' in \cite{GM98}.

\subsection{Survival probabilities}\label{survival4}

\noindent 
The case when the field $\xi$ assumes the values $-\infty$ and 0 only has 
a nice interpretation in terms of survival probabilities and is therefore 
of particular importance. The fundamental papers \cite{DV75} and \cite{DV79} 
by Donsker and Varadhan on the \mbox{Wiener} sausage contain apparently the 
first substantial annealed results on the asymptotics for the parabolic 
Anderson model. In the nineties, the thorough and deep work by Sznitman 
(see his monograph \cite{S98}), pushed the rigorous understanding of the 
quenched situation much further. 

\subsubsection{Brownian motion in a Poisson field of traps} 
We consider the continuous case, i.e., the version of \eqref{Anderson4} with 
$\Z^d$ replaced by $\R^d$ and the lattice \mbox{Laplacian} replaced by the 
usual \mbox{Laplace} operator. The field $\xi$ is given as follows. Let 
$(x_i)_{i\in I}$ be the points of a homogeneous Poisson point process in 
$\R^d$, and consider the union ${\mathcal O}$ of the balls $B_a(x_i)$ 
of radius $a$ around the Poisson points $x_i$.  We define a random potential
by putting 
\begin{equation}
\xi(x)=\begin{cases}0&\mbox{if }x\notin\Ocal,\\
-\infty&\mbox{if }x\in\Ocal.
\end{cases}
\end{equation}
(There are more general versions of this type of potentials, but for simplicity 
we keep with that.)  The set ${\mathcal O}$ receives the meaning of the set of 
`hard traps' or `obstacles'. Let 
$T_{\mathcal O}=\inf\{t>0\colon X(t)\in{\mathcal O}\}$ 
denote the entrance time into ${\mathcal O}$ for a \mbox{Brownian} motion 
$(X(t))_{t\in[0,\infty)}$. Then we have the \mbox{Feynman}-\mbox{Kac} 
representation 
$$
u(t,x)=\P_0\left(T_{\mathcal O}>t,\,X(t)\in dx\right)/dx, 
$$ 
i.e., $u(t,x)$ is equal to the sub-probability density of $X(t)$ on survival in 
the Poisson field of traps by time $t$ for \mbox{Brownian} motion starting 
from the origin. The total mass $U(t)=\P_0(T_{\mathcal O}>t)$ is the survival 
probability by time $t$. It is easily seen that the first moment of $U(t)$ 
coincides with a negative exponential moment of the volume of the \mbox{Wiener} 
sausage $\bigcup_{s\in[0,t]}B_a(X(s))$. 

\mbox{Donsker} and \mbox{Varadhan} analyzed the leading asymptotics of 
$\langle U(t)\rangle$ by using their large deviation principle for 
\mbox{Brownian} occupation time measures. The relevant islands have radius of 
order $\alpha(t)=t^{1/(d+2)}$. To handle the quenched asymptotics of $U(t)$, 
\mbox{Sznitman} developed a coarse-graining scheme for \mbox{Dirichlet} 
eigenvalues on random subsets of $\R^d$, the so-called method of enlargement of 
obstacles (MEO). The MEO replaces the eigenvalues in certain complicated 
subsets of $\R^d$ by those in coarse-grained subsets belonging to a 
{\em discrete\/} class of much smaller combinatorial complexity 
such that control is kept on the relevant properties of the eigenvalue. 

Qualitatively, the considered model falls into the class of bounded from above 
fields introduced in Section~\ref{distinc4} with $\gamma=0$. 

Related potentials critically rescaled with time have been studied in 
particular by \mbox{van den Berg} et al.\ \cite{BBH01} and by \mbox{Merkl} and 
\mbox{W\"uthrich} \cite{MW02}.

\subsubsection{Simple random walk among Bernoulli traps}
This is the discrete version of \mbox{Brownian} motion among \mbox{Poisson} 
traps. Consider the i.i.d.\ field $\xi=(\xi(x))_{x\in\Z^d}$ where $\xi(x)$ 
takes the values 0 or $-\infty$ only. Again, $u(t,x)$ is the survival 
probability of continuous-time random walk paths from $0$ to $x$ among the 
set of traps ${\mathcal O}=\{y\in\Z^d\colon \xi(y)=-\infty\}$. 

In their paper \cite{DV79}, Donsker and Varadhan also investigated the discrete 
case and described the logarithmic asymptotics of $\langle U(t)\rangle$ by 
proving and exploiting a large deviation principle for occupation times of 
random walks. Later \mbox{Bolthausen} \cite{B94} carried out a  deeper analysis 
of $\langle U(t)\rangle$ in the two-dimensional case using refined large 
deviation arguments. \mbox{Antal} \cite{Ant94}, \cite{Ant95} developed a 
discrete variant of the MEO and demonstrated its value by proving limit 
theorems for the survival probability $U(t)$ and its moments.

\subsection{General fields bounded from above}\label{distinc4}

\noindent
In \cite{BK01a} and \cite{BK01b}, a large class of potentials with 
$\esssup\xi(0)<\infty$ is considered. Assume for simplicity that 
$\esssup\xi(0)=0$ and that the tail of $\xi(0)$ at $0$  is given by
\begin{equation}\label{dist4}
\Prob\bigl(\xi(0)>-x\bigr)\approx\exp\left\{-D
x^{-\frac{\gamma}{1-\gamma}}\right\}, \qquad x\downarrow 0, 
\end{equation}
with $D>0$ and $\gamma\in [0,1)$ two parameters. The case $\gamma=0$ contains 
simple random walk among \mbox{Bernoulli} traps as a particular case. 
The cumulant generating function is roughly 
$H(t)\approx -\const t^\gamma$, and the annealed scale function is 
$\alpha(t)\approx t^\nu$ where $\nu=(1-\gamma)/(2+d-d\gamma)$. The power $\nu$ 
ranges from 0 to $1/(d+2)$ as the parameter $\gamma$ ranges from 1 to 0. 

It turns out in \cite{BK01a} that the rate function $I_R$ (see \eqref{LDPpot}) 
is given by 
\begin{equation}
I_R(\varphi)=\const \int_{Q_R}|\varphi(x)|^{-\frac \gamma{1-\gamma}}\,dx,
\end{equation}
where in the case $\gamma=0$ we interpret the integral as the \mbox{Lebesgue} 
measure of the support of $\varphi$. The characteristic variational formula 
for the annealed field shapes in \eqref{chideffirst} has been analyzed in great 
detail in the case $\gamma=0$. In particular it was shown that the minimizer is 
unique and has compact support, and it was characterized in terms of 
\mbox{Bessel} functions. However, in the general case $\gamma\in(0,1)$, 
an analysis of \eqref{chideffirst} has not yet been carried out.

\subsection{Gaussian fields and Poisson shot noise}\label{GauPo}

\noindent 
Two important particular cases in the {\em continuous\/} version of the 
parabolic Anderson model are considered in \cite{GK98} and \cite{GKM99} 
(see also \cite{CM95} for first rough results). The continuous version of 
\eqref{Anderson4} replaces $\Z^d$ by $\R^d$ and the discrete lattice 
\mbox{Laplacian} by the usual \mbox{Laplace} operator. Unlike in the discrete 
case, where any distribution on $\R$ may be used for the definition of an 
i.i.d.\ potential, in the continuous case it is not easy to find examples 
of fields that can be expressed in easily manageable terms. Since a certain 
degree of regularity of the potential is required, the condition of spatial  
independence must be dropped. 

In \cite{GK98} and \cite{GKM99}, two types of fields are considered: a 
\mbox{Gaussian} field $\xi$ whose covariance function $B$ has a parabolic 
shape around zero with $B(0)=\sigma^2>0$, and a so-called \mbox{Poisson} 
shot-noise field, which is defined as the superposition of copies of 
parabolic-shaped {\em positive\/} clouds around the points of a 
homogeneous \mbox{Poisson} point process in $\R^d$ (in contrast to the trap 
case of Section~\ref{survival4}). A certain (mild) assumption on the decay of 
the covariance function (respectively of the cloud) at infinity ensures 
sufficient independence between regions that are far apart.

Both fields easily develop very high peaks on small islands 
(the \mbox{Poisson} shot noise field is large where many \mbox{Poisson} 
points are close together). 
The annealed scale function is $\alpha(t)=t^{-1/4}$ for the \mbox{Gaussian} 
field and $\alpha(t)=t^{d/8}e^{-\sigma^2t/4}$ for the \mbox{Poisson} field 
\cite{GK98}.

\section{Results for the double-exponential case}\label{Results}

In this section, we formulate our results on the large-time asymptotics 
of the parabolic Anderson model in the particularly important case of a 
double-exponentially distributed random potential, see 
Section~\ref{doubleexp4}. We handle the annealed asymptotics of the total 
mass $U(t)$ in Section~\ref{Annres}, the quenched ones in 
Section~\ref{Quenres}, and the geometric picture of intermittency in 
Section~\ref{Interres}. The material of the first two subsections is
taken from \cite{GM98}, that of the last subsection from \cite{GKM04}.

\subsection{Annealed asymptotics}\label{Annres}

\noindent
As before, we assume that $\xi=(\xi(x))_{x\in\mathbb Z^d}$ is a field 
of i.i.d.\ random variables. We impose the following assumption on the 
cumulant generating function of $\xi(0)$ defined by (\ref{Hdef}). 
\vspace{2ex}

\noindent{\bf Assumption~(H).} {\em
The function $H(t)$ is finite for all $t>0$. There exists 
$\varrho\in[0,\infty]$ such that
$$
\lim_{t\to\infty} \frac{H(ct)-cH(t)}{t} = \varrho c\log c
\qquad \text{for all $c\in(0,1)$.}
$$
}
\vspace{2ex}

Note that the vicinity of the double-exponential distribution 
(\ref{doubleexp4}) corresponds to $\varrho\in(0,\infty)$. If $\varrho=\infty$, 
then the upper tail of the distribution of $\xi(0)$ is heavier than in 
the double exponential case, whereas for $\varrho=0$ it is thinner. 

Let $\mathcal P(\mathbb Z^d)$ denote the space of probability measures on 
$\mathbb Z^d$. We introduce the \mbox{Donsker}-\mbox{Varadhan} functional 
$S_d$ and the entropy functional $I_d$ on $\mathcal{P}(\mathbb Z^d)$ by
$$
S_d(\mu) = \sum_{\genfrac{}{}{0pt}{}{\{x,y\}\subset\mathbb Z^d}{|x-y|=1}}
\left( \sqrt{\mu(x)} - \sqrt{\mu(y)} \right)^2
\quad\text{and}\quad
I_d(\mu) = - \sum_{x\in\mathbb Z^d} \mu(x)\log \mu(x),
$$
respectively, and set
\begin{equation}\label{chid}
\chi_d= \inf_{\mu\in\mathcal P(\mathbb Z^d)}
\left[\kappa S_d(\mu) + \varrho I_d(\mu) \right],
\qquad \varrho\in[0,\infty].
\end{equation}

As before, let $U(t)$ denote the total mass of the solution $u(t,\cdot)$ to 
the PAM (\ref{Anderson4}). 

\begin{theorem}
Let Assumption~(H) be satisfied. Then, for any $p\in\mathbb N$, 
$$
\langle U(t)^p \rangle =
\exp\left\{ H(pt) - \chi_d pt + o(t) \right\}
\qquad \text{as $t\to\infty$.}
$$
\end{theorem}

It turns out that $\chi_d=2d\kappa\chi_0(\varrho/\kappa)$, where 
$\chi_0\colon[0,\infty)\to[0,1)$ is strictly increasing and concave,
$\chi_0(0)=0$, and $\chi_0(\varrho)\to 1$ as $\varrho\to\infty$. Moreover, 
for $\varrho\in(0,\infty)$, each minimizer $\mu$ of the variational problem 
(\ref{chid}) has the form 
$\mu=\const v^2$, where $v=v_1\otimes\dots\otimes v_d$
and each of the factor $v_1,\dots,v_d$ is a positive solution of the 
equation
$$
\kappa\Delta v + 2\varrho v\log v = 0 
\qquad \text{on $\mathbb Z$}
$$
with {\em minimal\/} $\ell^2$-norm. Uniqueness of $v$ modulo shifts 
holds for large $\varrho/\kappa$ but is open for small values of this 
quantity.  

Recall that $U(t)=\hat{u}(t,0)$, where 
$\hat{u}$ is the solution to (\ref{Anderson4}), but with non-localized 
initial datum $\1$ instead of $\delta_0$. 
A much deeper question is the computation of the asymptotics of the 
`correlation'
$$
c(t,x) = \frac{\langle\hat{u}(t,0)\hat{u}(t,x)\rangle}{\langle\hat{u}(t,0)^2
\rangle}
$$
of the spatially homogeneous solution $\hat{u}$ of the PAM. Assuming
additional regularity of the cumulant generating function $H$ and 
uniqueness modulo spatial shifts of the minimizer $\mu$ of the variational 
problem (\ref{chid}), it was shown in \cite{GH99} that  
$$
\lim_{t\to\infty} c(t,x) = 
\frac{\sum_z v(z)v(z+x)}{\sum_z v(z)^2}\,.
$$ 
This indicates that the second moment (considered as 
$\lim_{R\to\infty} |B_R|^{-1}\sum_{x\in B_R} \hat{u}^2(t,x)$) 
is generated by rare high peaks of the 
solution $\hat{u}(t,\cdot)$ {\em with shape\/} $v(\cdot)$.

\subsection{Quenched asymptotics}\label{Quenres}

Here we again consider i.i.d.\ random potentials $(\xi(x))_{x\in\mathbb Z^d}$ 
in the vicinity of the double-exponential distribution (\ref{doubleexp4}) 
but formulate our assumptions in a different manner. 

To be precise, let $F$ denote the distribution function of $\xi(0)$.
Provided that $F$ is continuous and $F(r)<\infty$ for all $r\in\R$
(i.e., $\xi$ is unbounded from above), we may introduce the non-decreasing function
\eq\label{phidef}
\varphi(r)=\log\frac 1{1-F(r)},\qquad r\in\R.
\end{equation}
Its left-continuous inverse $\psi$ is given by
\eq\label{psidef}
\psi(s)=\min\{r\in\R\colon \varphi(r)\geq s\},\qquad s>0.
\end{equation}
Note that $\psi$ is strictly increasing with $\varphi(\psi(s))=s$ 
for all $s>0$. The relevance of $\psi$
comes from the observation that $\xi$ has the same distribution as 
$\psi\circ\eta$, where
$\eta=(\eta(x))_{x\in\Z^d}$ is an i.i.d.\ field of exponentially 
distributed random variables with mean one.

We now formulate our main assumption.
\medskip

\noindent{\bf Assumption~(F). }{\em The distribution function 
$F$ is continuous, $F(r)<1$ for all $r\in\R$, and, in dimension $d=1$, 
$\int_{-\infty}^{-1}\log|r|\,F(\d r)<\infty$. There exists 
$\varrho\in(0,\infty]$ such that 
\begin{equation}\label{psiassa}
\lim_{s\to\infty} [\psi(cs)-\psi(s)]=\varrho\log c,\qquad 
c\in(0,1).
\end{equation}
If $\varrho=\infty$, then $\psi$ satisfies in addition 
\begin{equation}\label{psiassb}
\lim_{s\to\infty}[\psi(s+\log s)-\psi(s)]=0.
\end{equation}
}
\medskip

The crucial supposition  \eqref{psiassa} specifies that the upper tail of the 
distribution of $\xi(0)$ is close to the double-exponential distribution 
\eqref{doubleexp4} for $\varrho\in(0,\infty)$ and is heavier for 
$\varrho=\infty$. Assumption~\eqref{psiassb} excludes too heavy tails. 
Note that \eqref{psiassb} is fulfilled for Gaussian but not for exponential 
tails.

The reader easily checks that \eqref{psiassa} implies that 
$\psi(t)\sim\varrho\log t$ as $t\to\infty$. Let 
\begin{equation}\label{htdef}
h_t=\max_{x\in B_t}\xi(x),\qquad t>0,
\end{equation}
be the height of the potential $\xi$ in $B_t=[-t,t]^d\cap \Z^d$. It can be 
easily seen that, under Assumption~(F), almost surely,
\begin{equation}\label{hasy}
h_t=\psi(d\log t)+o(1)\qquad \mbox{as } t\to\infty.
\end{equation}
Let us remark that it is condition \eqref{psiassb} which ensures that 
the almost sure asymptotics of $h_t$
in \eqref{hasy} is non-random up to  order $o(1)$.

One of the main results in \cite{GM98}, Theorem~2.2, is the  second order 
asymptotics of the total mass $U(t)$ defined in \eqref{Udef}. 

\begin{theorem}
Under Assumption~(F), with probability one,
\begin{equation}\label{earlier}
\log U(t)=t\, [h_{t}-\widetilde{\chi}_d + o(1)]
\qquad \text{as $t\to\infty$.} 
\end{equation}
\end{theorem}

Here $0\le\widetilde{\chi}_d\le 2d\kappa$. 
An analytic description of $\widetilde{\chi}_d$ is as follows. 
Define $I\colon [-\infty,0]^{\Z^d}\to[0,\infty]$ by
\begin{equation}\label{Ldef}
I(V)=\begin{cases}
\sum_{x\in \Z^d} e^{V(x)/\varrho},&\text{if $\varrho\in(0,\infty)$,}\\
|\{x\in\Z^d\colon V(x)>-\infty\}|,&\text{if $\varrho=\infty$.}
\end{cases}
\end{equation}
\noindent
One should regard $I$ as {\em large deviation rate function\/} 
for the fields $\xi-h_t$ (recall (\ref{LDPpot}) and note that $\alpha(t)=1$ 
here). Indeed, if the distribution of $\xi$ is exactly given by 
\eqref{doubleexp}, then we have
$$
\Prob\left(\xi(\cdot)-h>V(\cdot) \mbox{ in }\Z^d\right)
=\exp\left\{-e^{h/\varrho}I(V)\right\}
$$
for any $V\colon\Z^d\to[-\infty,0]$ and any $h\in(0,\infty)$.
For $V\in[-\infty,0]^{\Z^d}$, let $\lambda(V)\in[-\infty,0]$ be
the top of the spectrum of the self-adjoint operator $\kappa\Delta+V$ 
in the domain $\{V>-\infty\}$ with zero boundary condition. In terms of the 
\mbox{Rayleigh}-{Ritz} formula,
\begin{equation}\label{RRform}
\lambda(V)=\sup_{f\in\ell^2(\Z^d)\colon\|f\|_2= 1}\bigl\langle
(\kappa\Delta+V)f,f\bigr\rangle,
\end{equation}
where $\langle\cdot,\cdot\rangle$ and $\|\cdot\|_2$ denote the inner product 
and the norm in $\ell^2(\Z^d)$, respectively. Then
\begin{equation}\label{chidef}
-\widetilde{\chi}_d=\sup\left\{\lambda(V)\colon V\in
[-\infty,0]^{\mathbb Z^d},\,I(V)\leq 1\right\}. 
\end{equation}
This variational problem is `dual' to the variational problem (\ref{chid}) 
and, in particular, $\widetilde{\chi}_d=\chi_d$.

\subsection{Geometry of intermittency}\label{Interres}

\noindent 
In this section we give a precise formulation of the geometric picture of 
intermittency, which was heuristically explained in Section~\ref{Interheu}. 

We keep the assumptions of the last subsection. For our deeper investigations, 
in addition to Assumption~(F), we introduce an assumption about the 
{\em optimal potential shape}.

\medskip
\noindent
{\bf Assumption~(M).} {\em Up to spatial shifts, the variational problem in 
\eqref{chidef} possesses a unique maximizer, which has a unique maximum.} 
\medskip

By $V_*$ we denote the unique maximizer of \eqref{chidef} which attains its 
unique maximum at the origin. We will call $V_*$ 
{\em optimal potential shape}. 
Assumption~(M) is satisfied at least for large $\varrho/\kappa$. 
This fact as well as further important properties of the 
variational problem \eqref{chidef} are stated in the next proposition.

\bep\label{phiprop}
(a) For any $\varrho\in(0,\infty]$, 
the supremum in \eqref{chidef} is attained.
\vspace{-1ex}
\begin{enumerate}
\item[(b)] If $\varrho/\kappa$ is sufficiently large, then the maximizer in 
\eqref{chidef} is unique modulo shifts and has a unique maximum.
\item[(c)] If Assumption~(M) is satisfied, then the optimal potential shape 
has the following properties.
\begin{enumerate}
\item[(i)] If $\varrho\in(0,\infty)$, then $V_*= 
f_*\otimes\cdots\otimes f_*$ for some $f_*\colon\Z\to(-\infty,0)$.
If $\varrho=\infty$, then $V_*$ is degenerate in the sense that 
$V_*(0)=0$ and $V_*(x)=-\infty$ for $x\not=0$. 
\item[(ii)] The operator $\kappa\Delta+V_*$ has a unique nonnegative 
eigenfunction
$w_*\in\ell^2(\Z^d)$ with $w_*(0)=1$ corresponding to the
eigenvalue $\lambda(V_*)$. Moreover, $w_*\in\ell^1(\Z^d)$. 
If $\varrho\in(0,\infty)$, then 
$w_*$ is positive on $\Z^d$, while $w_*=\delta_0$ for $\varrho=\infty$. 
\end{enumerate}
\end{enumerate}
\ep

We shall see that the main contribution to the total mass $U(t)$ comes from a 
neighborhood of the set of best local coincidences of $\xi-h_t$ with spatial 
shifts of $V_*$. These neighborhoods are widely separated from each other and 
hence not numerous. We may restrict ourselves further to those neighborhoods 
in which, in addition, $u(t,\cdot)$, properly normalized, is close to $w_*$. 

Denote by $B_R(y)=y+B_R$ the closed box of radius $R$ centered at $y\in\Z^d$ 
and write
\begin{equation}
B_R(A)=\bigcup_{y\in A}B_R(y)
\end{equation}
for the `$R$-box neighborhood' of a set $A\subset \Z^d$. In particular, 
$B_0(A)=A$.

For any $\eps> 0$ and any sufficiently large $\varrho\in(0,\infty]$, let
$r(\eps,\varrho)$ denote the smallest $r\in\N_0$ such that
\begin{equation}\label{rchoice}
\|w_*\|_2^2 \sum_{x\in\Z^d\setminus B_r}w_*(x)< \eps.
\end{equation}
Note that $r(\eps,\infty)=0$, due to the degeneracy of $w_\infty$.
Given $f\colon \mathbb Z^d\to\mathbb R$ and $R>0$, let 
$\|f\|_R = \sup_{x\in B_R}|f(x)|$. 

The main result of \cite{GKM04} is the following.

\begin{theorem}\label{main} Let the Assumptions~(F) and~(M) be satisfied. 
Then there exists a random time-dependent subset
$\Gamma^*=\Gamma^*_{t\log^2t}$ of $B_{t\log^2t}$ 
such that, almost surely, 
\begin{eqnarray}
(i)&&\label{mainresult}
\liminf_{t\to\infty}\frac 1{U(t)}\sum_{x\in B_{r(\eps,\varrho)}
(\Gamma^*)}u(t,x)\geq 1-\eps,\qquad \eps\in(0,1);\\
(ii)&&\label{largedistance}
|\Gamma^*|\leq t^{o(1)}\quad\mbox{and}
\quad \min_{y,\widetilde y\in\Gamma^*\colon y\not=
\widetilde y}|y-\widetilde y|\geq t^{1-o(1)}\qquad\mbox{as }t\to\infty;\\
(iii)&&\label{xiapprox}
\lim_{t\to\infty}\max_{y\in \Gamma^*}\,
\bigl\|\xi(y+\cdot)-h_t-V_*(\cdot)\bigr\|_R=0,
\qquad R>0;\\
(iv)&&\label{uapprox}
\lim_{t\to\infty}\max_{y\in \Gamma^*}\,\Bigl\|\frac{u(t,y+\cdot)}{u(t,y)}-
w_*(\cdot)\Bigr\|_R=0,
\qquad R>0.
\end{eqnarray}
\end{theorem}

Theorem~\ref{main} states that, up to an arbitrarily small relative error 
$\eps$, the islands with centers in $\Gamma^*$ and radius $r(\eps,\varrho)$ 
carry the whole mass of the solution $u(t,\cdot)$. Locally, in an arbitrarily 
fixed $R$-neighborhood of each of these centers, the shapes of the potential 
and  the normalized solution resemble $h_t+V_*$ and  $w_*$, 
respectively. The number of these islands increases at most as an arbitrarily 
small power of $t$ and their distance increases almost like $t$. Note that, 
for $\varrho=\infty$, the set $B_{r(\eps,\varrho)}(\Gamma^*)$ in 
\eqref{mainresult} is equal to $\Gamma^*$ and, hence, the islands consist of 
single lattice sites. 

It is an open problem under what assumptions on the potential
$\xi$ the number $|\Gamma^*|$ of relevant peaks stays 
bounded as $t\to\infty$; we have made no attempt
to choose $\Gamma^*$ as small as possible.

\section{Universality}\label{Universality}

\noindent 
In this section we explain that, under some mild regularity assumptions on 
the tails of $\xi(0)$ at its essential supremum, there are only four 
universality classes of asymptotic behaviors of the parabolic \mbox{Anderson} 
model. Three of them have already been analyzed in the literature: 
the double-exponential distribution  with $\varrho\in(0,\infty)$ respectively 
$\varrho=\infty$ \cite{GM98}, \cite{GH99}, \cite{GK98}, \cite{GKM99}, 
\cite{GKM04}) (see Sections~\ref{doubleexp} and \ref{Results}) and general 
fields bounded from above \cite{BK01a}, \cite{BK01b}, \cite{S98}, 
\cite{Ant94}, \cite{Ant95} (see Section~\ref{distinc4}). A fourth and new 
universality class is currently under investigation, see \cite{HKM04}. 
This class lies in the union of the boundary cases $\varrho=0$ of the 
double-exponential distribution (`almost bounded' fields) and $\gamma=1$ for 
the general bounded from above fields. Examples of distributions that fall 
into this class look somewhat odd, but it turns out that the optimal potential
shape and the optimal shape of the solution are perfectly parabolic, 
respectively \mbox{Gaussian}, which makes this class rather appealing. 
In particular, the appearing variational formulas can be easily solved 
explicitly and uniquely.

We now summarize \cite{HKM04}. Our basic assumption on the logarithmic moment 
generating function $H$ in \eqref{Hdef} is the following.

\medskip
\noindent 
{\bf Assumption~(\^H):} {\em There are a function 
$\widehat{H}\colon(0,\infty)\to\R$ and a continuous auxiliary function 
$\eta\colon(0,\infty)\to(0,\infty)$ such that} 
\begin{equation}\label{basic}
\lim_{t\uparrow\infty} \frac{H(ty)-yH(t)}{\eta(t)} = \widehat{H}(y)\not=0 
\qquad\text{for $y\not=1$.}
\end{equation}
\medskip

The function $\widehat H$ extracts the asymptotic scaling properties of the 
cumulant generating function $H$. In the language of the theory of regular 
functions, the assumption is that the logarithmic moment generating function 
$H$ is in the \mbox{de Haan} class, which does not leave many possibilities 
for $\widehat H$:

\begin{prop}\label{variation}
Suppose that Assumption~(\^H) holds. 
\begin{itemize}
\item[(i)] There is a $\gamma\ge 0$ such that 
$\lim_{t \uparrow \infty}\eta(yt)/\eta(t)= y^\gamma$ 
for any $y>0$, i.e., $\eta$ is regularly varying of index $\gamma$. 
In particular, $\eta(t)=t^{\gamma+o(1)}$ as $t\to\infty$. 
\item[(ii)] There exists a parameter $\rho >0$ such that, for every $y>0$,
\begin{itemize}
\item[(a)]  $\displaystyle\widehat H(y)=\rho \, \frac{y-y^\gamma}{1-\gamma}$ 
if $\gamma\not= 1$, 
\item[(b)] $\widehat H(y)=\rho  y \log y$ if $\gamma=1$. 
\end{itemize}
\end{itemize}\end{prop}

Our second regularity assumption is a mild supposition on the auxiliary 
function $\eta$. This assumption is necessary only in the case $\gamma=1$ 
(which will turn out to be the critical case).

\medskip
\noindent
{\bf Assumption~(K):} {\em The limit 
$\eta_*=\lim_{t \to \infty} \eta(t)/t\in[0,\infty]$ 
exists.}
\medskip

We now introduce a scale function $\alpha\colon(0,\infty)\to(0,\infty)$, by
\begin{equation}\label{alphadefannealed}
\frac{\eta\bigl(t\alpha(t)^{-d}\bigr)}{t\alpha(t)^{-d}}=\frac 1{\alpha(t)^2}\,.
\end{equation}
The function $\alpha(t)$ turns out to be the annealed scale function for the 
radius of the relevant islands in the parabolic \mbox{Anderson} model. 
We can easily say something about the asymptotics of $\alpha(t)$:

\begin{lemma}\label{alphaprop}
Suppose that Assumptions~(\^H) and~(K) hold.
If $\gamma\le 1$ and $\eta_*<\infty$, then there exists a unique solution 
$\alpha\colon(0,\infty)\to(0,\infty)$ to \eqref{alphadefannealed},
and it satisfies $\lim_{t\to\infty} t{\alpha(t)^{-d}}=\infty.$ Moreover, 
\begin{itemize}
\item[(i)] If $\gamma=1$ and $0<\eta_*<\infty$, then 
$\lim_{t\to\infty} \alpha(t)=1/\sqrt{\eta_*}\in(0,\infty)$.
\item[(ii)] If $\gamma=1$ and $\eta_*=0$, then $\alpha(t)=t^{\nu+o(1)}$ as 
$t\to\infty$, where $\nu=(1-\gamma)/(d+2-d\gamma)$. 
\end{itemize}
\end{lemma}

Now, under Assumptions~(\^H) and~(K), we can formulate a complete 
distinction of the PAM into four cases: 

\begin{itemize}
\item[(1)] $\eta_*=\infty$ (in particular, $\gamma\geq 1$).\\
This is the boundary case $\varrho=\infty$ of the double-exponential case. 
We have $\alpha(t)\to 0$ as $t\to\infty$, as is seen from 
\eqref{alphadefannealed}, i.e., the relevant islands consist of single lattice 
sites.
\item[(2)] $\eta_*\in(0,\infty)$ (in particular, $\gamma=1$).\\
This is the case of the double-exponential distribution in 
Section~\ref{doubleexp}. By rescaling, one can achieve that $\eta_*=1$. 
The parameter $\varrho$ of Proposition~\ref{variation}(ii)(b) is identical to 
the one in Assumption~(H) of Section~\ref{Annres}. 

\item[(3)] $\eta_*=0$ and $\gamma=1$.\\
This is the case of islands of slowly growing size, i.e., 
$\alpha(t)\to\infty$ as $t\to\infty$ slower than any power of $t$. This case 
comprises `almost bounded' and bounded from above potentials. This class is 
the subject of \cite{HKM04}, see also below. 
\item[(4)] $\gamma<1$ (in particular, $\eta_*=0$)\\
This is the case of islands of rapidly growing size, i.e., 
$\alpha(t)\to\infty$ as $t\to\infty$ at least as fast as some power of $t$. 
Here the potential $\xi$ is necessarily bounded from above. This case was 
treated in \cite{BK01a}; see Section~\ref{distinc4}.
\end{itemize}

Let us comment on the class (3), which appears to be new in the literature 
and is under 
investigation in \cite{HKM04}. One obtains examples of potentials (unbounded 
from above) that fall into this class by replacing $\varrho$ in 
\eqref{doubleexp4} by a sufficiently regular function $\varrho(r)$ that tends 
to $0$ as $r\to\infty$, and other examples (bounded from above) by replacing 
$\gamma$ in \eqref{dist4} by a sufficiently regular function $\gamma(x)$ 
tending to 1 as $x\downarrow 0$. According to Lemma~\ref{alphaprop}, the scale 
function $\alpha(t)$ defined in \eqref{alphadefannealed} tends to infinity, 
but is slowly varying. The annealed rate function for the rescaled potential 
shape, $I_R$, introduced in (\ref{LDPpot}) turns out to be
\begin{equation}
I_R(\varphi)=\const \int_{Q_R} e^{\varphi(x)/\varrho}\,d x.
\end{equation}
 The characteristic variational problem for the annealed potential shape in \eqref{chideffirst}, $\chi$, turns out to be uniquely minimized by a parabolic function $\varphi_*(x)={\rm const}-\varrho \|x\|_2^2$, and the principal eigenfunction $v_*$ of the operator $\kappa\Delta+\varphi_*$ is the Gaussian density $v_*(x)=\const e^{-\varrho \|x\|_2^2}$.

\section{Time-dependent random potentials}\label{Catsection}

\noindent
In this section we study the intermittent behavior of the parabolic 
\mbox{Anderson} model (PAM) with a space-time homogeneous ergodic
random potential $\xi$:
\begin{equation}\label{PAMt}
\begin{array}{rcll}
\displaystyle
\partial_t \,u(t,x)\!\!\! &=&\!\!\!\kappa \Delta u(t,x)
+\left[\xi(t,x)-\langle \xi(t,x)\rangle\right] 
u(t,x),\qquad &(t,x)\in(0,
\infty)\times {\mathbb Z}^d,\\
u(0,x)\!\!\!&=&\!\!\!1, &x\in{\mathbb Z}^d.
\end{array}
\end{equation}
Note that for time-dependent potentials the direct connection to the
spectral representation of the \mbox{Anderson} \mbox{Hamiltonian} 
(\ref{AndHam}) is lost. Our focus will be on the situation 
when the potential $\xi$ is given by a field $\{Y_k(t);k\in\mathbb N\}$ 
of independent random walks on $\mathbb Z^d$ with diffusion constant
$\varrho$ in \mbox{Poisson} equilibrium with density $\nu$:
\begin{equation}\label{catalyst}
\xi(t,x) = \gamma \sum_{k\in\mathbb N} \delta_{Y_k(t)}(x),
\end{equation}
where $\gamma$ denotes a positive coupling constant. Clearly
$$
\langle \xi(t,x) \rangle = \nu \gamma.
$$
This deterministic correction to the potential in (\ref{PAMt})
has been added for convenience to eliminate non-random terms. 

The form (\ref{catalyst}) of the potential is motivated by the following
particle model. Consider a system of two types of independent particles, 
$A$ and $B$, performing independent continuous-time simple random walks 
with diffusion constants $\kappa$ and $\varrho$ and \mbox{Poisson} 
initial distribution with densities $\nu$ and $1$, respectively. 
Assume that the $B$-particles split into two at a rate that is $\gamma$ times 
the number of $A$-particles present at the same location and die at rate 
$\nu\gamma$. Hence, the $A$- and $B$-particles may be regarded as 
catalysts and reactants in a simple catalytic reaction model. Then the 
(spatially homogeneous and ergodic) solution $u(t,x)$ of the PAM (\ref{PAMt}) 
is nothing but the average number 
of reactants at site $x$ at time $t$ given a realization of the catalytic 
dynamics (\ref{catalyst}). Such a particle model (with arbitrary death 
rate) has been considered by \mbox{Kesten} and \mbox{Sidoravicius} 
\cite{KS03}. We will come back to the results in \cite{KS03} at the end of 
this section. We further refer to the overview papers by \mbox{Dawson} and 
\mbox{Fleischmann} \cite{DF00} and by \mbox{Klenke} \cite{Kl00} for continuum 
models with singular catalysts in a measure-valued context where questions 
different from ours have been addressed.  

Our aim is to study the moment \mbox{Lyapunov} exponents 
$$
\lambda_p = \lim_{t\to\infty} \frac{1}{t} \log\,
            \langle u(t,0)^p \rangle
$$ 
as well as the quantities
$$
\lambda_p^* = \lim_{t\to\infty} \frac{1}{t} \log \log\, 
              \langle u(t,0)^p \rangle
$$
($p=1,2,\dots$) as functions of the model parameters. The phase diagram will  
turn out to be different in dimensions $d=1,2$, $d=3$, and $d\ge 4$. 

\begin{definition}\rm
a) For $p\in\mathbb N$, we will say that the PAM 
(\ref{PAMt})--(\ref{catalyst}) is
{\em strongly $p$-catalytic\/} if $\lambda_p^*>0$. Otherwise 
the PAM will be called {\em weakly $p$-catalytic}. 

b) For $p\in\mathbb N\setminus\{1\}$ and $\lambda_p<\infty$, we will say that 
the PAM is {\em $p$-intermittent\/} if 
$\lambda_p/p > \lambda_{p-1}/(p-1)$.
\end{definition}
We believe that strongly catalytic behavior is related to heavy tails of 
the \mbox{Poisson} distribution of the catalytic point process 
$\{Y_k(t); k\in\mathbb N\}$ and may occur if the main contribution to the 
$p$-th moment comes from 
realizations with a huge number of catalysts at the same lattice site.
Recall that $p$-intermittency means that, for large $t$, the $p$-th 
moment, considered as 
$\lim_{R\to\infty} |B_R|^{-1}\sum_{x\in B_R} u(t,x)^p$,   
is `generated' by high peaks of the solution $u(t,\cdot)$ located far 
from each other.

For potentials of the form
$$
\xi(t,x) = \gamma \dot{W}_x(t)
$$ 
with $(W_x(t))_{x\in\mathbb Z^d}$ being a field of i.i.d.\ 
(or correlated) \mbox{Brownian} motions and (\ref{PAMt}) understood as a 
system of \mbox{It\^o} equations, the moment \mbox{Lyapunov} 
exponents have been shown by \mbox{Carmona} and \mbox{Molchanov} \cite{CM94} 
to exhibit the following behavior.  
In dimensions $d=1,2$ there is $p$-intermittency for all 
$p\in\mathbb N\setminus\{1\}$ and all choices of the model parameters 
$\kappa$ and $\gamma$. If $d\ge 3$, then there exist critical points 
$0<c_2<c_3<\cdots$ such that $p$-intermittency holds if and 
only if $\kappa/\gamma<c_p$. For this model, the
asymptotics of the almost sure (`quenched') \mbox{Lyapunov} exponent 
as $\kappa\to 0$ has been investigated in \cite{CM94}, \cite{CMV96}, 
\cite{CKM01}, and \cite{CMS02}.

In the following we present the results for catalytic potentials of the  
form (\ref{catalyst}) obtained in \cite{GH04}. 

Our analysis of the moment \mbox{Lyapunov} exponents is based on the following 
probabilistic representation of the $p$-th moment which is easily derived 
from the \mbox{Feynman}-\mbox{Kac} formula for $u(t,0)$:
\begin{equation}\label{FKmoment}
\langle u(t,0)^p \rangle = \mathbb E_0^{(p)}
\exp\left\{ \nu\gamma \int_0^t \sum_{i=1}^p w(s,X_i(s))\,ds \right\},
\end{equation}
where $\mathbb E_0^{(p)}$ denotes expectation with respect to $p$ independent 
random walks $X_1,\dots,X_p$ on $\mathbb Z^d$ with generator 
$\kappa\Delta$ starting at the origin, and $w$ is the solution of the 
random initial value problem
\begin{equation*}
\begin{array}{rcll}
\displaystyle
\partial_t \,w(t,x)\!\!\! &=&\!\!\!\varrho \Delta w(t,x)
+\gamma\sum_{i=1}^p \delta_{X_i(t)}(x) \left(1+w(t,x)\right),
\qquad &(t,x)\in(0,
\infty)\times {\mathbb Z}^d,\\
w(0,x)\!\!\!&=&\!\!\!0, &x\in{\mathbb Z}^d.
\end{array}
\end{equation*}
One of the main difficulties in analyzing (\ref{FKmoment}) is related to 
the observation that $w(t,X_i(t))$ depends in a nontrivial way  
{\em on the whole past\/} $\{X_j(s); 0\le s \le t\}$, $j=1,\dots,p$,  
of our random walks (although our notation does not reflect this). 

For $r\ge 0$, let $\mu(r)$ denote the upper boundary of the spectrum of the 
operator $\Delta+r\delta_0$ in $\ell^2(\mathbb Z^d)$. It is well-known 
that, in dimensions $d=1,2$, $\mu(r)>0$ for all $r$ and, 
in dimensions $d\ge 3$, $\mu(r)=0$ for $0\le r\le r_d$ and $\mu(r)>0$ for  
$r>r_d$, where
$$
r_d = 1/G_d(0)
$$
and $G_d$ denotes the \mbox{Green} function associated with the discrete 
\mbox{Laplacian} $\Delta$. 

\begin{theorem}\label{GdH1}
For any choice of the parameters of the PAM (\ref{PAMt})--(\ref{catalyst}), 
the limit $\lambda_p^*$ exists, and
$$
\lambda_p^* = \varrho\, \mu(p\gamma/\varrho),
\qquad p\in\mathbb N.
$$
\end{theorem}
Hence, for $d=1,2$, the PAM (\ref{PAMt})--(\ref{catalyst}) is always 
strongly $p$-catalytic, whereas 
for $d\ge 3$ this is true only if $p\gamma/\varrho$ exceeds the critical 
threshold $r_d$. Note that $\lambda_p^*$ does not depend on $\kappa$ nor 
on $\nu$.

We next study the behavior of the moment \mbox{Lyapunov} exponents 
$\lambda_p=\lambda_p(\kappa)$ as a function of the diffusion constant 
$\kappa$ in the weakly catalytic regime $0<p\gamma/\varrho<r_d$ for 
$d\ge 3$. We will mainly describe their behavior for small and for 
large values of $\kappa$. 

\begin{theorem}\label{GdH2}
Let $d\ge 3$, $p\in\mathbb N$, and $0<p\gamma/\varrho<r_d$. Then the 
limit $\lambda_p$ exists and is finite for all $\kappa$ and $\nu$. 
Moreover, $\kappa\mapsto\lambda_p(\kappa)$ is strictly decreasing and 
convex on $[0,\infty)$ and satisfies
$$
\lim_{\kappa\downarrow 0} \frac{\lambda_p(\kappa)}{p}
= \frac{\lambda_p(0)}{p}
= \nu\gamma \frac{p\gamma/\varrho}{r_d-p\gamma/\varrho}.
$$ 
\end{theorem}
Hence, in any dimension $d\ge 3$, the PAM (\ref{PAMt})--(\ref{catalyst}) is 
$p$-intermittent in the weakly catalytic regime for small values of the 
diffusion constant $\kappa$. 

To formulate the behavior of the moment \mbox{Lyapunov} exponents for 
$\kappa\to\infty$, we introduce the variational expression $\mathcal P$ 
for the polaron problem analyzed in \cite{L77}, \cite{DV83}, and 
\cite{BDS93}:
$$
\mathcal P 
= \sup_{\genfrac{}{}{0pt}{}{f\in C_c^\infty(\mathbb R^3)}{\|f\|_2=1}} 
\left[ \left\| (-\Delta)^{-1/2} f^2 \right\|_2^2 
- \left\| \nabla f \right\|_2^2 \right].
$$
\vspace{1ex}

\begin{theorem}\label{GdH3}
Let $d\ge 3$, $p\in\mathbb N$, and $0<p\gamma/\varrho<r_d$. Then
$$
\lim_{\kappa\to\infty} \kappa\frac{\lambda_p(\kappa)}{p} = 
\begin{cases}
\frac{\nu\gamma^2}{r_3}+\sqrt{p}\sqrt{\frac{\nu\gamma^2}{\varrho}} \mathcal P,
&\text{if $d=3$,}\\
\frac{\nu\gamma^2}{r_d},
&\text{if $d\ge 4$.}
\end{cases}
$$
\end{theorem}
In other words, in dimensions $d=3$ for $p$ in the weakly catalytic regime, 
the PAM (\ref{PAMt})--(\ref{catalyst}) is $p$-intermittent also for 
large values of $\kappa$. We conjecture 
intermittent behavior for all $\kappa$. In dimensions $d\ge 4$, 
the leading term in the asymptotics of $\lambda_p(\kappa)/p$ as 
$\kappa\to\infty$ is the same for all $p$ in the weakly catalytic regime. 
We in fact conjecture that there is even no intermittency in high dimensions 
for large $\kappa$. 

Let us remark that in \cite{KS03} \mbox{Kesten} and \mbox{Sidoravicius} 
obtained the following results for the above mentioned catalytic particle 
model with arbitrary death rate $\delta$ instead of $\nu\gamma$. If $d=1,2$, 
then for all choices of the parameters, the average 
number of $B$-particles per site tends to infinity faster than exponential.
If $d\ge 3$, $\gamma$ small enough and $\delta$ large enough, then the 
average number of $B$-particles per site tends to zero exponentially fast.
The first of these two results is covered by Theorem~\ref{GdH1}, 
since the death rate $\delta$  does not affect $\lambda_1^*$. 
The second result is covered by Theorem~\ref{GdH2}, 
since $0<\lambda_1<\infty$ in the weakly catalytic regime and the death 
rate $\delta$ shifts $\lambda_1$ by $\gamma\nu-\delta$. 
\mbox{Kesten} and \mbox{Sidoravicius} further show that in all dimensions
for $\gamma$ large enough, conditioned on the evolution of the $A$-particles,
there is a phase transition. Namely, for small $\delta$ the $B$-particles
locally survive, while for large $\delta$ they become locally extinct. 
In \cite{GH04} there are no results for the quenched situation. The 
analysis in \cite{KS03} does not lead to an identification of \mbox{Lyapunov} 
exponents, but is more robust under an adaption of the model than the 
above analysis based on the \mbox{Feynman}-\mbox{Kac} representation.

\end{document}